\documentclass[12pt]{article}
\usepackage{amsmath, dsfont, amssymb, amscd, caption,  color}
\usepackage{graphicx,subfigure}
\usepackage[margin=1in]{geometry}

\usepackage{times,soul}
\usepackage[square,sort,comma,numbers]{natbib}

\usepackage{setspace}
\usepackage{tikz}
\usepackage{url}
\usepackage{mathbbol}
\usepackage{tikz-cd}

\usepackage[pagewise]{lineno}

\usepackage{amsfonts}
\usepackage{fancybox}
\usepackage{algorithmic}
\usepackage{algorithm}
\usepackage{comment}
\usepackage{color}
\usepackage[colorlinks,citecolor=blue,urlcolor=blue]{hyperref}

\usepackage{paralist}
\usepackage{tikz,enumerate}
\usepackage{multirow}
\usepackage[normalem]{ulem}

\usepackage{hyperref}
\usepackage{cleveref}

\newcommand{\bpm}{\begin{pmatrix}}
\newcommand{\epm}{\end{pmatrix}}

\graphicspath{{images/}}

\newtheorem{theorem}{Theorem}
\newtheorem{lemma}{Lemma}
\newtheorem{prop}{Proposition}
\newtheorem{corollary}{Corollary}

\newtheorem{remark}{Remark}
\newtheorem{example}{Example}
\newtheorem{definition}{Definition}

 \numberwithin{dummy}{section}

\newenvironment{proof}[1][. ]{{\bf Proof#1}}{\hfill$\square$\vskip\baselineskip}

\begin{document}

\begin{center}
    
{\Large\bf{Cohomology of Quaternionic Foliations and Orbifolds}}

{\normalsize \vspace{1cm} }

{\normalsize \vspace{0.5cm} Rouzbeh Mohseni \\[0pt]
Faculty of Mathematics and Computer Science, University of Łódź \\[0pt]
ul. Banacha 22, 90-238 Łódź, Poland \\[0pt]
e-mail: rouzbeh.mohseni@wmii.uni.lodz.pl \\[0pt]
https://orcid.org/0000-0003-3773-6268}
 
{\normalsize \vspace{0.5cm} Robert A. Wolak \\[0pt]
Institute of Mathematics, Jagiellonian University \\[0pt]
ul. \L ojasiewicza 6, 30-348 Krak\'{o}w, Poland \\[0pt]
e-mail: robert.wolak@uj.edu.pl\\[0pt]
https://orcid.org/0000-0002-7587-4998}

\end{center}

\date{}

\begin{abstract}
Starting with a concise review of quaternionic geometry and quaternionic K{\"a}hler manifolds, we define a transversely quaternionic K{\"a}hler foliation. Then we formulate and prove the foliated versions of the now classical results of V.Y. Kraines and A. Fujiki on the cohomology of quaternionic K{\"a}hler manifolds. Finally, as any orbifold can be realized as the leaf space of a suitably defined Riemannian foliation, we reformulate our results for quaternionic orbifolds.


\end{abstract}

\section{Introduction}\label{sec1}

In recent decades, several geometric theories have been successfully developed and implemented in modern physics, thus strengthening the relationship between these two sciences. String theory is one of them, and the appearance of extra dimensions in it, though initially striking and somewhat embarrassing, later led to new approaches and tools that, in a sense, were created to justify these ``undetectable" extra dimensions. The manifolds that can be candidates for the geometry of these extra dimensions are compact Einstein manifolds that have positive scalar curvature. Depending on the theory one is working on, whether superstring theory (i.e., $10$ dimensions) or M-theory (i.e., $11$ dimensions), subtracting the $4$ dimensions of spacetime, we are left with $6$ or $7$ dimensions. Geometries including Calabi-Yau, hyperk\"ahler, quaternionic K\"ahler, $G_2$ and Spin(7) manifolds are automatically Einstein and have many applications in physics, e.g., supersymmetric sigma models and supergravity, see \cite{al,bg3} and references therein. Looking closely, it appears that all of these geometries appear in Berger's list of holonomy groups \cite{ber} with K\"ahler geometry as the main one missing. Although a K\"ahler manifold is not automatically Einstein, there exist K\"ahler-Einstein manifolds, with Calabi-Yau as one of the most important examples. An important class of odd-dimensional geometries that also appears in physics is Sasakian-Einstein and $3$-Sasakian geometry. Let $(M, g)$ be a Sasakian-Einstein manifold, then its cone $(C(M), \tilde{g}=ds^2+s^2g)$ is K\"ahler and Ricci-flat, the converse is also true. As for a $3$-Sasakian manifold, the cone is hyperk\"ahler. These odd-dimensional geometries not only have a relation to K\"ahler and hyperk\"ahler geometries through their cones but also can appear as a bundle over a quaternionic K\"ahler manifold (orbifold), cf. \cite{bg1, bg2, bg3}. Moreover, Sasakian and $3$-Sasakian manifolds also admit transversely K\"ahler and quaternionic K\"ahler foliations respectively and are related to all of the geometries we mentioned earlier, see \cite{bg3}.

V.Y. Kraines in \cite{kr} gave an analog of the Hodge decomposition theorem for a quaternionic manifold (back then there was no distinction between what is now known as quaternionic K\"ahler and the Ricci-flat case, i.e. hyperk\"ahler). Next, using some results of S.S. Chern of \cite{ch}, she demonstrated inequalities on Betti numbers. Later A. Fujiki \cite{fu} formulated analogs of the Hodge and Lefschetz decompositions theorems for the cohomology of some special manifolds, in particular quaternionic K{\"a}hler manifolds. In this paper, we investigate these results for the case of a foliated quaternionic K{\"a}hler structure on a Riemannian manifold. In Section \ref{sec2}, for the convenience of the reader, we gather some definitions and results deemed necessary to understand the theory developed in the final two sections. First we discuss quaternionic K{\"a}hler manifolds and quaternionic K{\"a}hler analogs of the Hodge star operator $*$ as well as the operators $L$ and $\Lambda$. In Subsection \ref{sec4}, we provide a basic introduction to foliated and transverse geometric structures associated to a foliation. Then we formulate the definition of a transversely quaternionic K{\"a}hler foliation. Next, we discuss Riemannian submersions over quaternionic K\"ahler manifolds, since they provide us with examples of the type of foliation we defined earlier. The final subsection is dedicated to twistor spaces. Twistor spaces over a manifold have proven to be important in the study of the properties of the base manifold. In particular, in the case of quaternionic K{\"a}hler manifolds the twistor space was studied by S. Salamon \cite{sa} and independently by L. B{\'e}rard Bergery \cite{br,be}. Salamon shows how using cohomology groups of the base manifold, some of the characteristic classes of the twistor space can be computed; hence, the study of these spaces can be fruitful in both directions. Having this in mind,  we define the transverse quaternionic twistor space $Z\mathcal{F}$ on a foliated manifold $(M^{p+4n},\mathcal{F})$ of codimension $q=4n$. 

The main results of the paper are presented in Section \ref{sec8}. First, for a transversely quaternionic K\"ahler foliation, the analogs of the operators $L$ and $\Lambda$ are introduced, using which we define basic effective forms. Then we prove a decomposition theorem for basic forms (Theorem \ref{thm3}).   The basic cohomology for transversely quaternionic K{\"a}hler foliations and the foliated counterparts of Kraines' results are the main results of this section. The estimates for the basic Betti numbers are formulated in Theorem \ref{thm5}. The injectivity of the mapping $L$ and the decomposition of basic cohomology groups are studied in Theorem \ref{th5}. The important result of A. Haefliger et al. \cite{gi} ensures that for any orbifold, one can find a Riemannian foliated manifold with compact leaves such that its leaf space is the original orbifold. The characteristic foliation of a 3-Sasakian is a nice example of transversely quaternionic K\"ahler foliations with compact leaves, and hence the leaf space has the structure of an orbifold. This fact motivated us to reformulate our results for orbifolds in the final section. Lastly, we have prepared three appendices that are organized as follows; In Appendix \ref{App1}, we recall fundamental notions of the theory of quaternionic manifolds, namely the $4$-form $\Omega$ defined by Kraines, which is used to define two operators $L$ and $\Lambda$ on the space of all forms.  In Appendix \ref{sec5}, we discuss basic forms and basic Hodge theory on foliated Riemannian manifolds. Finally in Appendix \ref{App3}, we recall the basics of the theory of foliated G-structures and their foliated connections.

\section{Preliminaries}\label{sec2}

For the convenience of the reader, in this section, we recall some basic
definitions and constructions that will be used in the following sections. We would like to point out that in each subsection we have a few new definitions and results that we shall use in the subsequent sections.

\subsection{Quaternionic Manifolds}\label{sec3}

The geometry of complex manifolds, in particular K\"ahler manifolds, is one of the most studied subjects in differential geometry. A K\"ahler manifold can be defined in many ways, one of which is through the reduction of the holonomy group to $U(n)$, with $n$ being the complex dimension of the manifold. In \cite{ber} Berger listed the possible holonomy groups for simply connected, irreducible, nonsymmetric Riemannian manifolds. At the time, $SO(n)$ and $U(n)$ were the only holonomy groups that were well-studied, which correspond respectively to the orientable and K\"ahler manifolds. Later, manifolds with holonomy group $Sp(n).Sp(1)$ were studied by Bonan \cite{bo} and Kraines \cite{kr}. Riemannian manifolds $(M, g)$ of real dimension 4n whose holonomy group can be reduced to $Sp(n).Sp(1)$ are called quaternionic K\"ahler. In dimension $4 (n=1)$ as $Sp(1).Sp(1) = SO(4)$, this condition is basically equivalent to the manifold being orientable and Riemannian; Due to this fact, we may assume $n \geq 2$. However, it should be noted that the case of $4$-dimensional manifolds is quite important in physics and the correct analog of the quaternionic K\"ahler condition for this case is for the manifold to be Einstein and self-dual \cite{be}. We will discuss the importance of this particular case later on. \newline

\noindent Quaternionic K{\"a}hler manifolds can also be characterized
in terms of local endomorphisms of the tangent bundle, cf., e.g.,
\cite[Proposition 14.36]{be}:

\begin{prop}
A Riemannian manifold $(M,g)$ is quaternionic K{\"a}hler if and only if there exists a covering of $M$ by open sets $U_i$ and, for each $i$, two almost complex structures $I$ and $J$ on $U_i$ such that
\begin{enumerate}[(i)]
    \item $g$ is Hermitian for $I$ and $J$ on $U_i$,
    \item $IJ=-JI$,
    \item the Levi-Civita derivatives of $I$ and $J$ are linear combinations of $I$, $J$ and $K=IJ$,
    \item for any $x\in U_i\cap U_j$ the vector space of endomorphisms of $T_xM$ generated by $I$, $J$ and $K$ is the same for $i$ and $j$.
\end{enumerate}

\end{prop}

\noindent In fact, condition (iv) states that these local endomorphisms $I$, $J$ and $K$ defined on each open subset $U_i$ generate a global subbundle commonly denoted by $Q$ of $End(TM)$, which is parallel with respect to the induced action of the Levi-Civita connection. The almost complex structures $I, J$ and $K$ are transformed by $SO(3)$ on their respective domains of existence.

\noindent The proposition itself is a consequence of the fact that the subgroup $Sp(n).Sp(1)$ of $SO(4n)$ can be characterized as the group of orientation preserving linear isometries which leave invariant the $3$-dimensional subspace of endomorphisms of $\mathbb{R}^{4n}$ generated by right multiplication by $i$, $j$ and $k$ when $\mathbb{R}^{4n}$ is identified with $\mathbb{H}^n$. $Sp(n)$ acts on the left by $(n,n)$-quaternion matrices and $Sp(1)$ acts on the right by multiplication by quaternions of norm $1$. Unfortunately, the endomorphisms $I$, $J$ and $K$ cannot be globally defined. Denote by $\omega_I, \omega_J$ and $\omega_K$ the corresponding K\"ahler $2$-forms of $I$, $J$ and $K$, respectively, then the following $4$-form

\begin{equation}\label{omegaQK}
     \Omega=\Omega_I\wedge \Omega_I+\Omega_J\wedge \Omega_J+\Omega_K\wedge \Omega_K.
\end{equation}
is a globally defined and is sometimes referred to as the Kraines $4$-form. Next we define operators $*$, $L$, and $\Lambda$ on the space $A^*(M)$ of differential forms on the manifold $M$:

\begin{flalign*}
    &*\colon A^k(M)\rightarrow A^{4n-k}(M),\\
   & L\colon A^k(M)\rightarrow A^{k+4}(M);\,\,\, L(\alpha)=\Omega\wedge \alpha,\\
&\Lambda\colon A^k(M)\rightarrow A^{k-4}(M);\,\,\, \Lambda(\alpha)=*(\Omega\wedge *\alpha).
\end{flalign*}
A differential form $\alpha$ is called \emph{effective} if $\Lambda\alpha=0$.

Since at any point $x$  of a quaternionic K\"ahler manifold $M$, the tangent space $T_xM$ can be identified with $\mathbb{H}^n$,  one can define a global closed $4$-form $\Omega$ of maximal rank by pulling back the form $\Omega\in (\mathbb{H}^n)^\prime$ defined in (\ref{omegaH}), and then show that it is the same form as that of  (\ref{omegaQK}). Moreover, Theorem \ref{Thm1} applied point by point, allows us to formulate and prove the following decomposition theorem for differentiable forms on $M$, cf. \cite[Theorem 3.5]{kr}.

\begin{theorem}
Let $M$ be a $4n$-dimensional quaternionic K{\"a}hler manifold and $\omega$ a differential form on $M$ of degree $p\leq n+1$. Then
\begin{align*}
    \omega=\sum^{[p/4]}_{i=0}L^i\omega_e^{p-4i},
\end{align*}
where $\omega^k_e$ is an effective $k$-form.
\end{theorem}

\subsection{Foliations}\label{sec4}

Let $\mathcal{F}$ be a foliation on Riemannian $m$-manifold $(M,g)$ of codimension $q$ and of leaves of dimension $p$. Then $\mathcal{F}$ is defined by a cocycle $\mathcal{U}=\{U_i,f_i,g_{ij}\}_{i,j \in I}$ modeled on a $q$-manifold $N_0$ such that

\begin{enumerate}
    \item $\{U_i\}_{i\in I}$ is an open covering of $M$,
    \item $f_i:U_i\rightarrow N_0$ are submersions with connected fibers,
    \item $g_{ij}:N_0\rightarrow N_0$ are local diffeomorphisms of $N_0$ with $f_i=g_{ij}f_j$ on $U_i\cap U_j$.
\end{enumerate}

\noindent The connected components of the trace of any leaf of $\mathcal{F}$ on $U_i$ consist of
the fibers of $f_i$. The open subsets $N_i = f_i(U_i) \subset N_0$ form a $q$-manifold
$N_{\mathcal{U}} = \bigsqcup N_i$, which can be considered as a transverse manifold of the
foliation $\mathcal{F}$. The pseudogroup $\mathcal{H}_\mathcal{U}$ of local diffeomorphisms of $N_{\mathcal{U}}$ generated
by $g_{ij}$ is called the holonomy pseudogroup of the foliated manifold
$(M,\mathcal{F})$ defined by the cocycle $\mathcal{U}$. Two different cocycles can define the
same foliation, then we have two different transverse manifolds and two
holonomy pseudogroups. In fact, these two holonomy pseudogroups are
equivalent in the sense of Haefliger, cf. \cite{ha}.

According to Haefliger, cf. \cite{ha}, a transverse property of a foliated
manifold is a property of foliations which is shared by any two foliations
with equivalent holonomy pseudogroups. For example, being Riemannian,
transversely symplectic, transversely almost-complex, transversely
Kähler, etc., is a transverse property. A Riemannian foliation, i.e.,
admitting a bundle-like metric, is defined by a cocycle $\mathcal{U}$ modeled on
a Riemannian manifold whose local submersions are Riemannian submersions.
Then the associated transverse manifold $N_{\mathcal{U}}$ is Riemannian
and the associated holonomy pseudogroup $\mathcal{H}_\mathcal{U}$ is a pseudogroup of local
isometries. Any foliation defined by a cocycle $\mathcal{V}$ whose holonomy pseudogroup
$\mathcal{H}_\mathcal{V}$ is equivalent to $\mathcal{H}_\mathcal{U}$ is also Riemannian, as an equivalence
of pseudogroups transports the Riemannian metric from $N_\mathcal{U}$ to $N_\mathcal{V}$ and
ensures that the pseudogroup $\mathcal{H}_\mathcal{V}$ is a pseudogroup of local isometries
of the transported metric. This metric can be lifted to a bundle-like
metric (not unique) on the other foliated manifold making the second
foliation Riemannian. The same procedure can be applied to any
geometrical structure, for the discussion of this general procedure see
\cite{wo-t}.

\begin{definition}
A foliation $\mathcal{F}$ is transversely quaternionic K{\"a}hler if it
is defined by a cocycle \break $\mathcal{U}=\{U_i,f_i,g_{ij}\}_{i,j \in I}$ modeled on a quaternionic K{\"a}hler manifold $(N_0, g_0, Q_0)$ and the local diffeomorphisms $g_{ij}$ are local
automorphisms of the quaternionic Kähler structure of $(N_0, g_0, Q_0)$, i.e.,
the $g_{ij}$ are local isometries and the induced mappings $\Tilde{g}_{ij}$ of $End(TN_0)$
preserve the subbundle $Q_0$ of rank $3$.
\end{definition}

\noindent In the language of foliated structures this condition can be formulated as follows, cf. \cite{pv}. Let $N(M,\mathcal{F})=TM/T\mathcal{F}\,$ be the normal bundle of the foliation $\mathcal{F}$. The vector bundle $End(N(M,\mathcal{F}))$ admits the natural foliation $\mathcal{F}_{End}$ of dimension $p$ which is defined by a cocycle $\mathcal{F}_{End}=\{V_i,\Tilde{f}_i,\Tilde{g}_{ij}\}_{i,j\in I}$ modeled on $End(TN_0)$ where $\Tilde{f}(A)=df\circ A\circ (df|_{N(M,\mathcal{F})})^{-1}$. With this in mind we can define a foliated quaternionic K{\"a}hler structure.

\begin{definition}
A foliated quaternionic Kähler structure on a foliated
Riemannian manifold $(M,\mathcal{F}, g)$ is given by the following data:
\begin{enumerate}
    \item $g$ is a foliated Riemannian metric in $N(M,\mathcal{F})$;
    \item a $3$-dimensional foliated subbundle $Q$ of $End(N(M,\mathcal{F}))$ which is locally spanned by $3$ foliated almost complex structures;
   \item the metric $g$ is Hermitian with respect to these local almost complex
structures;
   \item the subbundle $Q$ is parallel with respect to the foliated Levi-Civita
connection of $g$.
\end{enumerate}
\end{definition}
Therefore, a foliated quaternionic Kähler structure on a foliated Riemannian
manifold $(M,\mathcal{F})$ will be denoted by $(M, \mathcal{F}, g, Q)$.
Let $g$ be a foliated Riemannian metric on $N(M,\mathcal{F})$ and $\Bar{g}$ the corresponding holonomy invariant metric on the transverse manifold $N$.
At each point $x \in U_i$ there exist $3$ foliated almost complex structures
$I_x$, $J_x$, and $K_x$ on an open neighbourhood $U_x$ which project to $3$ almost complex structures $\Bar{I}_x$, $\Bar{J}_x$, and $\Bar{K}_x$ on a neighborhood of $f_i(x) \in N_i$. Then on $U_x$ we define the
$2$-forms
\begin{align*}
    \Omega_I(u, v) = g(Iu, v),\Omega_J (u, v) = g(Ju, v), \text{and}\,\,\, \Omega_K(u, v) = g(Ku, v),
\end{align*}
\noindent
where $u, v \in N(M,\mathcal{F})$.
Using the same argument as in \cite{kr} one can show that the $4$-form $\Omega$
\begin{align*}
    \Omega=\Omega_I\wedge\Omega_I+\Omega_J\wedge\Omega_J+\Omega_K\wedge\Omega_K
\end{align*}
is well-defined, i.e., it is independent of the choice of the structures
$I$, $J$, and $K$. As these structures were foliated the form $\Omega$ is basic.
Moreover, $\Omega$ is closed and of maximal rank and parallel with respect to
the foliated Levi-Civita connection.
In exactly the same way using the transverse metric $\Bar{g}$ we define the $4$-form $\Bar{\Omega}$ of the transverse manifold $N$. The $4$-form $\Bar{\Omega}$ is closed and of maximal rank. It is also holonomy invariant as the Riemannian metric $\Bar{g}$ is. The forms $\Omega$ and $\Bar{\Omega}$ correspond to each other under the correspondence between foliated and transverse objects, cf. \cite{wo-t}.

\subsection{Riemannian submersions over quaternionic K\"ahler manifolds}\label{sub2.3}

\begin{definition}
    Let $M^{2n+1}$ be a smooth manifold together with a structure $(\phi, \xi, \eta, g)$, where $\phi$ is a $(1,1)$-tensor field of rank $2n$, $\xi$ is a vector field and $\eta$ is $1$-form such that
    \begin{align}
        \phi^2=-1+\eta\otimes\xi,\,\,\,\,\, \eta(\xi)=1,
    \end{align}
    and $g$ is a compatible Riemannian metric, i.e., it satisfies $g(\phi X, \phi Y)=g(X, Y)-\eta(X)\eta(Y),\break\forall X,Y\in\Gamma(TM)$. Then $(M^{2n+1}, \phi, \xi, \eta, g)$ is called an \emph{almost contact metric manifold}.
\end{definition}
\noindent
An almost contact metric structure is called:
\begin{itemize}
\setlength{\itemindent}{-.1in}
\item \emph{contact metric structure}, if it satisfies $d\eta=2g(X, \phi Y)$,
    \item \emph{normal}, if it satisfies $N_\phi+2d\eta\otimes\xi=0$, where $N_\phi$ is the Nijenhuis tensor of $\phi$,
    \item \emph{Sasakian}, if it is a normal contact metric structure.
\end{itemize}

\begin{remark}\label{remark1}
    It is well-known that on a Sasakian manifold, the Reeb vector field $\xi$ defines a $1$-dim foliation $\mathcal{F}_\xi$ called the characteristic or Reeb foliation. If this foliation has only compact leaves (i.e., quasi-regular case), its leaf space can be a manifold (i.e., regular case) or an orbifold. It is common to refer to such cases as quasi-regular (the leaf space is an orbifold) or regular (the leaf space is a manifold) Sasakian manifolds. If some leaves are not compact, the term irregular is used. We will discuss the regular case in more detail shortly and postpone the quasi-regular case to Section \ref{sec9}. For more discussion, see \cite{bg2}.
\end{remark}

\begin{definition}
    A smooth manifold $M^{4n+3}$ equipped with three almost contact structures $(\phi_i, \xi_i, \eta_i),\break  i=1,2,3$, is called an \emph{almost $3$-contact manifold} if 
    \begin{align}
        \phi_k=\phi_i\phi_j-\eta_j\otimes\xi_i,\,\,\,\,\, \xi_k=\phi_i\xi_j,\,\,\,\,\eta_k=\eta_i\circ\phi_j,
    \end{align}
    holds for any cyclic permutation of $\{1,2,3\}$.
\end{definition}
\noindent
It is well-known that on such a manifold there always exists a Riemannian metric $g$ that satisfies $g(\phi_i X,\phi_i Y)=g(X, Y)-\eta_i(X)\eta_i(Y),\,i=1, 2, 3$, and $(\phi_i, \xi, \eta_i, g)$ is then called an \emph{almost \break$3$-contact metric structure}, cf. \cite{kuo}. The tangent bundle decomposes as $TM=\mathcal{V}\oplus \mathcal{H}$, where $\mathcal{V}$ is spanned by $\xi_i$, hence $\mathcal{H}$ is of rank $4n$ and these distributions are called vertical and horizontal, respectively.

\noindent One can investigate the normality of each of these three structures and, if all are normal, then the $3$-structure is said to be \emph{hypernormal}. Traditionally, the notion of \emph{$3$-Sasakian} manifolds was as follows: Three distinct Sasakian structures $(\phi_i, \xi_i, \eta_i)$ on a smooth manifold $M^{4n+3}$ with a Riemannian metric that is compatible with all three and they satisfy
\begin{align}\label{3v}
    \xi_k=\frac{1}{2}[\xi_i, \xi_j],\,\,\, g(\xi_i, \xi_j)=\delta_{ij}
\end{align}
for a cyclic permutation of $\{1, 2, 3\}$. However, Kashiwada \cite{ka} showed that on an almost $3$-contact manifold, if each of the three structures is a contact metric structure, then the manifold is necessarily $3$-Sasakian. The main reference for the theory of $3$-Sasakian manifolds is Ch. Boyer's and K. Galicki's book \cite{bg2} or their paper \cite{bg1}.\newline

\noindent
The quaternionic K\"ahler structure is closely related to the $3$-Sasakian structure, and it was the subject of studies in a series of papers by Ishihara and Konishi. It is obvious that the vertical distribution $\mathcal{V}$ spanned by the three vector fields $\xi_i$ on a $3$-Sasakian manifold $\tilde{M}$ is actually integrable due to (\ref{3v}), assuming that it is also regular, $\tilde{M}$ is a foliated manifold with foliation denoted by $\mathcal{F}$ and the space of leaves $M:=\tilde{M}/\mathcal{F}$ is a smooth manifold of dimension $4n$. The Riemannian submersion $\pi\colon\tilde{M}\rightarrow M$ and several differential objects on these two manifolds was the topic of their joint work in \cite{Is1}. Then, in \cite{Is2}, Ishihara demonstrated that $M^{4n}$ is in fact a quaternionic K\"ahler manifold. Later, Sakamoto \cite{sk}, proved that on every quaternionic K\"ahler manifold there exists an $SO(3)$-principal bundle $P$, and Konishi \cite{ko} constructed on $P$ a $3$-Sasakian structure, which is canonically associated to the quaternionic K\"ahler structure on the base. The characteristic folitaion of a $3$-Sasakian manifold is transversely quaternionic K\"ahler \cite[Theorem 2.8]{bg2} and its twistor space is a trivial bundle as a foliated bundle since on the normal bundle of the characteristic foliation we have 3 linearly independent foliated almost-complex structures.\newline

\noindent Recently, a new class of almost $3$-contact metric manifolds has been introduced by Agricola and Dileo \cite{Ag20}, which generalizes that of $3$-Sasakian manifolds and is called \emph{$3$-$(\alpha,\delta)$-Sasaki manifolds}.

\begin{definition}
    A $3$-$(\alpha,\delta)$-Sasaki manifold is an almost $3$-contact metric manifold $(M, \phi_i, \xi_i, \eta_i, g)$ such that
    \begin{align}
        d\eta_i=2\alpha\Phi_i+2(\alpha-\delta)\eta_j\wedge\eta_k,
    \end{align}
    where $\Phi_i(X, Y):=g(X,\phi_i Y)$, $\alpha\neq 0,\,\,\delta$ are real constants and indices are true for any cyclic permutaion of $\{1,2,3\}$.  
\end{definition}

Later in a series of papers, the properties of such manifolds and examples were studied \cite{Ag21,Ag23}, however, for our purposes here, we only mention a few. Notice that for $\alpha=\delta=1$, we have a $3$-Sasakian manifold.

\begin{theorem}\cite[Theorem 2.2.1]{Ag20}
Any $3$-$(\alpha,\delta)$-Sasaki manifold is hypernormal.
\end{theorem}

\noindent
Let $M$ be a $3$-$(\alpha,\delta)$-Sasaki manifold, then we have a locally defined Riemannian submersion $\pi:(M,g)\rightarrow (N,g_N)$, cf. \cite[Proposition 2.1.1]{Ag21} and on the base a quaternionic K\"ahler structure defined by
\begin{align}
    \check{\phi}_i=\pi_*\circ\phi_i\circ s_*,
\end{align}
where $s:N\rightarrow M$ is an arbitrary section of $\pi$, cf. \cite[Theorem 2.2.1]{Ag21}. The scalar curvature of the base manifold $N$ is $16n(n+2)\alpha\delta$, cf. \cite[Theorem 2.2.2]{Ag21}.
\begin{corollary}\cite[Corollary 2.3.1]{Ag20}
Let $(M,\phi_i, \xi_i, \eta_i, g)$ be a $3$-$(\alpha, \delta)$-Sasaki manifold. Then for every smooth vector field $X$ on $M$ and cyclic permutation of $\{1, 2, 3\}$,
\begin{equation}
    \nabla^{LC}_X\xi=-\alpha\phi_iX+(\alpha-\delta)[\eta_j(X)\xi_k-\eta_k(X)\xi_j],\,\,\,
    \nabla^{LC}_{\xi_i}\xi_i=0,\,\,\,\, \nabla^{LC}_{\xi_i}\xi_j=-\nabla^{LC}_{\xi_j}\xi_i=\alpha\xi_k,
    \end{equation}
hence each $\xi_i$ is a Killing vector field and the vertical distribution $\mathcal{V}$ is integrable with totally geodesic leaves.
    
\end{corollary}

\noindent Since the vertical distribution $\mathcal{V}$ is intergrable, it induces a folitaion on the manifold $M$ and as $\mathcal{V}$ is spanned by Killing vector fields, the metric $g$ is bundle-like \cite{bg1}. Therefore, the manifold $(M,\phi_i, \xi_i, \eta_i, g)$ is an interesting example of a Riemannian manifold with transversely quaternionic K\"ahler foliation.

\subsection{Twistor spaces}\label{sec7}

 Let $(M^{4n},g)$ be a quaternionic K{\"a}hler manifold and following notation of Salamon \cite{sa}, let $\bf{H}$ be the quaternionic line bundle associated to the representation of $Sp(1)$ on $\mathbb{C}^2$. Salamon showed that on such a quaternionic K{\"a}hler manifold $(M,g)$ there exists a twistor fibration $q:Z\rightarrow M$. $Z$ is a complex manifold of complex dimension $2n+1$ which can be viewed as an $S^2$-bundle generated by anticommuting almost complex structures $I,J,K$ or as $Z=\mathbb{P}(\bf{H})$, i.e. the projectivization of the bundle $\bf{H}$ and $Z$ can be refered to as the twistor space of the quaternionic K{\"a}hler manifold $M$. It is worth noting that $\bf{H}$ is not always globally well-defined but $Z=\mathbb{P}(\bf{H})$ is. Twistor space of quaternionic K{\"ahler} manifolds and in general quaternionic manifolds can be used to study some of the properties of the base manifold. In particular Salamon showed how some of the characteristic classes of the twistor space $Z$ can be computed using cohomology groups of the manifold $M$.\newline

\noindent Let $\mathcal{F}$ be a Riemannian foliation of codimension $4q$ on the manifold $M^{p+4q}$ and on $(M, \mathcal{F})$ we have a foliated quaternionic K{\"a}hler structure $Q$ with local bases $(I_1,I_2,I_3)$ on the normal bundle $N(M,\mathcal{F})=TM/T\mathcal{F}$.
We define the transverse quaternionic twistor space of $\mathcal{F}$ by
\begin{align*}
 Z\mathcal{F}=\{J\in Q,\,\, J=\alpha_1I_1+\alpha_2I_2+\alpha_3I_3,\,\, \alpha_1^2+\alpha_2^2+\alpha_3^2=1\} 
\end{align*}
i.e. $Z\mathcal{F}$ is the sphere bundle associated with the foliated vector bundle $Q$, and $Q$ carries a Riemannian structure such that it makes $\{I_1,I_2,I_3\}$ an orthonormal basis. 

\noindent In our earlier work \cite{rm}, we discussed the twistor spaces of foliated manifolds and the corresponding transverse manifold and how these two are related. Now as we discussed in Subsection \ref{sec4} the quaternionic K{\"a}hler foliation is defined by a cocycle modeled on a quaternionic  K{\"a}hler manifold $N_0$. Therefore one can expect that since the twistor space $Z(N_0)$ of  $N_0$ is a useful tool in studying its geometric properties, using the relation between $Z(N_0)$ and $Z\mathcal{F}$ it would be possible to investigate properties of the foliated manifold $(M,\mathcal{F})$.\newline

\section{Basic cohomology of transversely quaternionic K{\"a}hler foliations}\label{sec8}

In this section, we focus mainly on Riemannian manifolds with a transversely quaternionic K\"ahler foliation and present several results regarding the decomposition of basic forms on them; we advise the reader to consult Appendix B in case of unfamiliarity with the theory of basic cohomology.\newline

\noindent Using the 4-form $\Omega$, we define $L$ and $\Lambda$ operators on the complex of basic forms $A^*(M,\mathcal{F})$:
\begin{flalign*}
   & L\colon A^k(M,\mathcal{F})\rightarrow A^{k+4}(M,\mathcal{F});\,\,\, L(\alpha)=\Omega\wedge \alpha,\\
    &\Lambda\colon A^k(M,\mathcal{F})\rightarrow A^{k-4}(M,\mathcal{F});\,\,\, \Lambda(\alpha)=\Bar{*}(\Omega\wedge \Bar{*}\alpha),
\end{flalign*}
\noindent Basic forms that are annihilated by $\Lambda$ are called effective.

\noindent On a compact manifold with a taut foliation, one can define scalar
products $\langle .\,,.\,\rangle$ and $\langle .\,,.\,\rangle_b$ on $A^k(M)$ and $A^k(M,\mathcal{F})$, respectively, as
\begin{enumerate}
    \item $\langle \omega,\omega^\prime\rangle=\int_M *(\omega\wedge*\omega^\prime)=\int_M\omega\wedge*\omega^\prime,$
    \item  $\langle \omega,\omega^\prime\rangle_b=\int_M \Bar{*}(\omega\wedge\Bar{*}\omega^\prime)=\int_M\omega\wedge\Bar{*}\omega^\prime\wedge\chi_{\mathcal{F}}$.
\end{enumerate}
Using this scalar product we have for any $\omega \in A^k(M,\mathcal{F})$ and $\omega^\prime \in
A^{k+4}(M,\mathcal{F})$
\begin{align*}
    \langle L\omega,\omega^\prime\rangle_b=\langle\omega,\Lambda\omega^\prime\rangle_b.
\end{align*}
As the proof of the decomposition theorem, cf. \cite[Theorem 3.5]{kr}, was
based on the pointwise application of the decomposition theorem for
forms on the quaternions, the same argument is valid in the case of
basic forms.

\begin{theorem}\label{thm3}
Let $(M, g, \mathcal{F})$ be a $(4n+p)$-dimensional Riemannian foliated
manifold whose \break $p$-dimensional foliation $\mathcal{F}$ is transversely quaternionic
K{\"a}hler. Let $\omega$ be a basic differential form on $(M,\mathcal{F})$ of degree $p\leq n + 1$.
Then
\begin{align*}
    \omega=\sum^{[p/4]}_{i=0}L^i\omega^{p-4i}_e,
\end{align*}
where $\omega^k_e$ is an effective basic $k$-form.

\begin{proof}\normalfont
The operator $L$ is an isomorphism into and $\Lambda$ as the adjoint of $L$ is onto. The proof can be done using induction on $p$. It is easy to see that by definition of $\Lambda$ we have $\bigwedge^p=\bigwedge^p_e$ for $p=0,1,2,3$. Now assume that for $i<p$ the theorem holds we shall prove it for $i=p$. We want to show that the space $\bigwedge^p_e$ is the orthogonal complement of the subspace $L\bigwedge^{p-4}(\mathbb{H}^n)^\prime$ in $\bigwedge^p(\mathbb{H}^n)^\prime$. Let $\omega\in\bigwedge^p_e$ and $L\omega^\prime\in L\bigwedge^{p-4}(\mathbb{H}^n)^\prime$ then since $\Lambda \omega=0$ we have 
\begin{align*}
    (\omega,L\omega^\prime)=(\Lambda\omega,\omega^\prime)=0,
\end{align*}
therefore orthogonality is proved. Now let $\omega\in \bigwedge^p(\mathbb{H}^n)^\prime$ and $(\omega,L\omega^\prime)=0, \forall \omega^\prime\in\bigwedge^{p-4}(\mathbb{H}^n)^\prime$ then
\begin{align*}
    (\Lambda \omega, \omega^\prime)=0,
\end{align*}
hence we have $\Lambda \omega=0$ and by induction the proof is complete. 
\end{proof}
\end{theorem}

In what follows, we shall use the theory of foliated $ G$-structures. For a reader unfamiliar with the topic, we summarized the basics of the theory in Appendix C. Let $(M,\mathcal{F})$ be a compact Riemannian foliated manifold. Assume
that
\begin{enumerate}[1)]
    \item its foliated normal bundle $(N(M,\mathcal{F}),\mathcal{F}_N)$ admits a reduction to
a connected subgroup $G$ of $O(q)$,
    \item the corresponding foliated $G$-reduction $B(M,G;\mathcal{F})$ of the foliated frame bundle \break $L(M,\mathcal{F}_L)$ admits a foliated connection without torsion. This condition is equivalent to the vanishing of the structure
tensor of $B(M, G;\mathcal{F})$, cf. \cite{wo-b}.
\end{enumerate}

These two foliated conditions correspond to the following transverse
ones:
\begin{enumerate}[1')]
    \item the tangent bundle of the transverse manifold $N_\mathcal{U}$ admits a reduction
to a connected subgroup $G$ of $O(q)$,
    \item the corresponding $G$-reduction $B(N_\mathcal{U},G)$ of the frame bundle
$L(N_\mathcal{U})$ is holonomy invariant, i.e., $\mathcal{H}_\mathcal{U}$-invariant, and admits a connection without torsion. This condition is equivalent to the vanishing
of the structure tensor of the $G$-structure. Since a $G$-connection without
torsion is unique, it is $\mathcal{H}_\mathcal{U}$-invariant, cf. \cite{fuji, wo-b}.
\end{enumerate}

The fiber bundle $ \bigwedge^k N_x(M,\mathcal{F})^*$ can be understood as the associated bundle of $L(M, \mathcal{F}_L)$ with the standard fiber $ \bigwedge^k(R^{q*})$.The space of sections of this bundle we denote by $A^k(N).$  Since the normal frame bundle $L(M,\mathcal{F}_L)$ is foliated, the foliation $\mathcal{F}_L$ induces a foliation $ \mathcal{F}_L^k$ of the fiber bundle $ \bigwedge^k N_x(M,\mathcal{F})^*$. The space of basic $k$-forms $A^k(M,\mathcal{F})$  is a subspace of $A^k(N)$. If the normal frame bundle $L(M,\mathcal{F}_L)$ admits a foliated $G$-reduction $B(M, G;\mathcal{F})$, the bundle $ \bigwedge^k N_x(M,\mathcal{F})^*$ can be understood as the associated fiber bundle of  $B(M, G; \mathcal{F})$ with the standard fiber $\bigwedge^k(R^{q*})$. The natural induced foliations coincide. 

Let $W \subset \bigwedge^k(R^{q*})$ be an invariant subspace of $\bigwedge^k(R^{q*})$ under the standard action of $G$. There is the standard scalar product on $\bigwedge^k(R^{q*})$ for which the induced action of $G$ is isometric. The associated fiber bundle 
$\mathcal W$ of  $B(M, G; \mathcal{F})$ with the standard fiber $W$ can be understood as a foliated vector subbundle of the foliated vector bundle $ (\bigwedge^k N_x(M,\mathcal{F})^*,  \mathcal{F}_L^k).$ Therefore a k-differential form $\alpha$ which corresponds to a section of  $\mathcal W$ is said to be of type $W$. The space of these $\mathcal W$-valued sections we denote also by $\mathcal{W}$. The projection $P_W\colon A^k(N) \rightarrow \mathcal{W}$ sends basic forms
into basic forms as the operation is done point by point. Next, we show that the result of S.S. Chern in \cite{ch} can be reformulated for the basic Laplacian $\Delta_b$.


\begin{prop}[cf. \cite{ch}]

Let $W \subset \bigwedge^k(R^{q*})$ be an invariant subspace of $\bigwedge^k(R^{q*})$ under the standard action of $G$, $P_W$ be the projection $P_W\colon A^k(M,\mathcal{F}) \rightarrow \mathcal{W}$ and $\Delta_b$ be the basic Laplacian, then 
\begin{align*}
    P_W\Delta_b=\Delta_b P_W.
\end{align*}
Moreover, let $W_1, \ldots,W_s$ be irreducible invariant subspaces of $\bigwedge^k(R^{q*})$
for the action of the group $G$. Then if $\alpha$ is a harmonic basic $k$-form, the $k$-forms
$P_{W_1}\alpha, \ldots, P_{W_s}\alpha$ are basic and  harmonic. Moreover, if $\alpha$ is a basic $k$-form of type
$W$ so is the form $\Delta_b\alpha$.

\begin{proof}\normalfont
Our proof follows Chern's original proof, cf. \cite{ch}. The unique torsionless connection in the foliated $G$-structure is foliated, which in particular means that the tangent bundle to the foliation $\mathcal{F}_G$ is a subbundle of the horizontal bundle, cf. \cite{wo-b}. One of the key points of Chern's proof is the identification of forms on the base manifold with horizontal forms on the total space $B$ of the $G$-bundle. The previous remark about foliations ensures that basic forms on $(M,\mathcal{F})$ correspond to horizontal basic forms on $B(M, G; \mathcal{F})$. Having that in mind Chern's proof can be easily adapted to the foliated case. For the reader's convenience, we recall the main steps of the proof. \newline


\noindent Since the canonical torsionless connection of $B(M, G; \mathcal{F})$ is foliated the structure equation for the canonical $1$-form $\theta$ can be written as
\begin{align}\label{dw}
    d\theta^i=-\sum_{\rho,a} a^i_{\rho a}\theta^{a}\wedge \pi^\rho
\end{align}
where $\theta=(\theta^1,\ldots,\theta^q)$ is a local representation of the canonical $1$-form of the foliated bundle $B(M, G; \mathcal{F})$ and $\pi^i$ is the set of linearly independent left-invariant 1-forms of $G$ lifted to the total space of the $G$-bundle. The forms $\theta^i$ can be chosen to be orthogonal.

\noindent The foliated curvature form of the foliated connection can be written as
\begin{align*}
    \Omega^i=\frac{1}{2} \sum_{j,a} R^i_{ja}\theta^j\wedge \theta^a,\,\,\, R^i_{ja}=-R^i_{aj}.
\end{align*}
 Let us next introduce new tensors which shall be of use later on
\begin{align}
    S^i_{jal}=\sum_\rho a^i_{\rho j}R^\rho_{al}.
\end{align}

\noindent Since $G$ is a subgroup of $O(q)$, there exists a metric and using it we can lower or raise the indices, so we shall from now on use only subscripts and we have the following relations
\begin{align}
    S_{ijal}=-S_{jial}=-S_{ijla},\\
    S_{ijal}=S_{alij}.
\end{align}

\noindent Let $\alpha$ be a basic $k$-form which can be written as 
\begin{align}\label{alpha}
    \alpha=\frac{1}{k!}\sum_{i_1 ,\ldots ,i_k}\alpha_{i_1\ldots i_k}\theta^{i_1}\wedge\ldots\wedge\theta^{i_k}.
\end{align}
where $(\theta^1,\ldots,\theta^q)$ is the local basis and $\alpha_{i_1\ldots i_k}$ is anti-symmetric in any two of its indices. In order to compute the basic Laplacian $\Delta_b\alpha$, initially we need to compute $d\alpha$ and to do this, first using (\ref{dw}) and differentiating the terms in (\ref{alpha}) we get the relation
\begin{align}\label{dw2}
    d\alpha_{i_1\ldots i_k}+\sum^k_{l=1}\sum_\rho \alpha_{i_1\ldots i_{l-1}ji_{l+1}\ldots i_k}a^j_{\rho i_l}\pi^\rho=\sum_m \alpha_{i_1\ldots i_k|m}\theta^m.
\end{align}
the notation used for the right hand side is originally used by Chern and almost resembles the notation usually used for the derivative of the coefficients in an ordinary exterior derivative. Now we can write the exterior derivative as
\begin{align*}
    d\alpha=\frac{1}{k!}\sum_{i_1 ,\ldots ,i_k,j}\alpha_{i_1\ldots i_k|j}\theta^j\wedge\theta^{i_1}\wedge\ldots\wedge\theta^{i_k}.
\end{align*}

\noindent Next, to have anti-symmetric coefficients we write it as
\begin{align*}
    d\alpha=\frac{(-1)^k}{(k+1)!}\sum_{i_1 ,\ldots ,i_{k+1}}(\alpha_{i_1\ldots i_k|i_{k+1}}-\alpha_{i_{k+1}i_2\ldots i_k|i_{1}}-\cdots-\alpha_{i_1\ldots i_{k-1}i_{k+1}|i_{k}})\theta^{i_1}\wedge\ldots\wedge\theta^{i_{k+1}}.
\end{align*}

 \noindent We also have
\begin{align*}
    *\alpha=\frac{1}{k!(q-k)!}\sum_{i_1 ,\ldots ,i_q}\epsilon_{i_1\ldots i_k i_{k+1}\ldots i_q}\alpha_{i_1\ldots i_{k}}\,\,\theta^{i_{k+1}}\wedge\ldots\wedge\theta^{i_{q}}.
\end{align*}
where $\epsilon_{i_1\ldots i_q}$ is equal to +1(-1) if $i_1,\ldots , i_q$ form an even(odd) permutation of $1, \ldots , q$ and is otherwise equal to zero.

\noindent Taking into account that on our manifold we have taut foliation $\mathcal{F}$, i.e. a metric can be chosen without changing the transverse geometry such that the leaves are minimal submanifolds and hence $\kappa=0$, we can write the adjoint $\delta_b$ of $d$
\begin{align*}
    \delta\alpha=\frac{(-1)^k}{(k-1)!}\sum_{i_1 ,\ldots ,i_{k-1},j}\alpha_{i_1\ldots i_{k-1}j|j}\,\,\theta^{i_1}\wedge\ldots\wedge\theta^{i_{k-1}}.
\end{align*}

\noindent Further differentiating (\ref{dw2}) with some extra steps taken gives us
\begin{align*}
    \alpha_{i_{1}\ldots i_{k}|m|r}-\alpha_{i_1\ldots i_{k}|r|m}=\sum^k_{l=1}\sum_j \alpha_{i_1\ldots i_{l-1}ji_{l+1}\ldots i_k}S_{ji_l rm}
\end{align*}
and this formula helps us to shorten the expression for the basic Laplacian given as follows
\begin{align*}
    \Delta_b\alpha=-\frac{1}{k!}\sum_{i_1 ,\ldots ,i_k,j}\alpha_{i_1\ldots i_k|j|j}\,\,\theta^{i_1}\wedge\ldots\wedge\theta^{i_k}-\frac{1}{(k-1)!}\sum_{i_1 ,\ldots ,i_k,r,j}\alpha_{i_1\ldots i_{k-1}r}S_{rji_kj}\,\,\theta^{i_1}\wedge\ldots\wedge\theta^{i_k}\\+\frac{1}{(k-2)!}\sum_{i_1 ,\ldots ,i_k,r,j}\alpha_{i_1\ldots i_{k-2}rj}S_{ri_{k-1}ji_k}\,\,\theta^{i_1}\wedge\ldots\wedge\theta^{i_k}.
\end{align*}

\noindent Now we are facing an issue where the coefficients in this formula are not anti-symmetric. In order to solve this issue we need to introduce the following quantities
\begin{align}\label{S}
    S(i_1\ldots i_k,j_1\ldots j_k; r_1\ldots r_k, m_1\ldots m_k)= 
    \sum\epsilon(i_1\ldots i_k;n_1\ldots n_{k-1}g)\\\nonumber \times\epsilon(j_1\ldots j_k;n_1\ldots n_{k-1}h)\epsilon(r_1\ldots r_k;l_1\ldots l_{k-1}u)\epsilon(m_1\ldots m_k;l_1\ldots l_{k-1}v)S_{ghuv}
\end{align}
where $\epsilon(i_1\ldots i_k;j_1\ldots j_{k})$ is $+1$ or $-1$ if $j_1\ldots j_{k}$ is an even or odd permutation of $i_1\ldots i_{k}$ respectively and is zero otherwise, moreover all the indices run from $1$ to $q$ and the summation is over all the repeated ones. For the sake of brevity, we will use the notation $S((i) (j); (r) (m))$ instead of (\ref{S}), using which we define yet another quantity as follows
\begin{align}\label{S2}
    S((i) (r))=S(i_1\ldots i_{k},r_1\ldots r_{k})=\frac{1}{k!}\sum_{j_1 ,\ldots ,j_{k}} S(i_1\ldots i_k,j_1\ldots j_k; r_1\ldots r_k, j_1\ldots j_k),
\end{align}
which is anti-symmetric in the indices belonging to each one of the sets $i_1 ,\ldots ,i_{k}$ and $r_1 ,\ldots ,r_{k}$ and $ S((i) (r))= S((r) (i))$. Using (\ref{S2}) we are able to rewrite the expression for the basic Laplacian 
\begin{flalign}
    \Delta_b\alpha&=-\frac{1}{k!}\sum_{i_1 ,\ldots ,i_k,j}\alpha_{i_1\ldots i_k|j|j}\,\,\theta_{i_1}\wedge\ldots\wedge\theta_{i_k}\\\nonumber&-\frac{1}{(k!(q-k)!)^2}\sum_{i_1 ,\ldots ,i_k,r_1 ,\ldots ,r_k}\alpha_{i_1\ldots i_{k}}S(i_1 ,\ldots , i_k, r_1 ,\ldots ,r_k)\,\,\theta_{r_1}\wedge\ldots\wedge\theta_{r_k}
\end{flalign}
Let $\theta^{i_1}\wedge\ldots\wedge \theta^{i_k}$ be the basis of $\bigwedge^k(R^{q*})$ where $1\leq i_1<\ldots<i_k\leq q$. Suppose $W_1$ is an invariant subspace of $\bigwedge^k(R^{q*})$ under the standard action of $G$ and $W_2$ be its orthogonal space, then $W_2$ is also invariant under this action. There exist base vectors $\gamma_{1}, \ldots ,\gamma_{K}\subset\bigwedge^k(R^{q*})$, $K=\begin{pmatrix} q \\ k\end{pmatrix}$ which are related to $\theta^{i_1}\wedge\ldots\wedge\theta^{i_k}$ via an orthogonal transformation
\begin{align}\label{gam}
    \theta^{i_1}\wedge\ldots\wedge\theta^{i_k}=\sum^K_{\mu=1}g_{i_1\ldots i_k,\mu}\gamma_\mu,
\end{align}
such that $\gamma_{1},\ldots,\gamma_{\nu}$ and $\gamma_{\nu+1},\ldots,\gamma_{K}$ span the subspaces $W_1$ and $W_2$ respectively. From (\ref{gam}) one can obtain
\begin{align*}
    \gamma_\mu=\frac{1}{k!}\sum_{i_1,\ldots ,i_{k}}g_{i_1\ldots i_k,\mu}\,\,\theta^{i_1}\wedge\ldots\wedge\theta^{i_k}.
\end{align*}

\noindent Using (\ref{dw}) we compute the exterior derivative 
\begin{align*}
    d\gamma_\mu=\frac{1}{(k-1)!}\sum_{i_1,\ldots ,i_k,m,\rho}g_{i_1\ldots i_{k-1},m,\mu}a_{m\rho i_k}\,\,\pi^\rho\wedge\theta^{i_1}\wedge\ldots\wedge\theta^{i_k}\\=\frac{1}{(k-1)!}\sum_{i_1,\ldots ,i_{k-1},m,A,\lambda,\rho}g_{i_1\ldots i_{k-1},m,\mu}g_{i_1\ldots i_{k-1},A,\lambda}a_{m\rho A}\,\,\pi^\rho\wedge\gamma^\lambda,
\end{align*}
here $1\leq\mu\leq\nu$. The condition of $W_1$ being invariant under the action of $G$ can be obtained as
\begin{align}\label{cc}
    \sum_{i_1,\ldots,i_{k-1},m,A}g_{i_1\ldots i_{k-1}m,\mu}g_{i_1\ldots i_{k-1}A,\zeta}S_{mAjr},\,\, \nu+1\leq \zeta\leq K.
\end{align}

\noindent Now suppose that $\alpha$ is a basic $k$-form of type $W_1$ and write it as 
\begin{align*}
    \alpha=\sum_\mu \alpha_\mu\gamma_\mu=\frac{1}{k!}\sum_{i_1,\ldots , i_{k}, \mu}\alpha_\mu g_{i_1\ldots i_{k},\mu}\,\,\theta^{i_1}\wedge\ldots \wedge\theta^{i_k}.
\end{align*}
and the basic Laplacian is the following
\begin{align*}
    \Delta_b\alpha=-\sum\alpha_{\mu|j|j}\gamma_\mu-\frac{1}{(k!(k-1)!)^2}\sum\alpha_\mu g_{i_1\ldots i_{k},\mu}g_{r_1\ldots r_{k},\lambda}S((i) (r))\gamma_\lambda.
\end{align*}

\noindent Therefore in order for $\Delta_b\alpha$ to be of type $W_1$ the second term needs to vanish for $\nu+1\leq\lambda\leq K$, which follows from the properties of $S((i) (r))$ and (\ref{cc}).
\end{proof}
\end{prop}

Let $\{e_i\}_{i=1,\ldots, q}$ be a local orthonormal frame of the normal bundle $Q$ and $\{e^i\}_{i=1,\ldots, q}$ be its dual coframe. Let $\alpha\in A^k(M,\mathcal{F})$ be a basic $k$-form, the exterior derivative preserves the basic forms and the restriction to basic forms denoted by $d_b=d|_{A^*(M,\mathcal{F})}$ is well-defined. It is possible to write $d_b$ and its adjoint $\delta_b$ for $\kappa=0$ as follows
\begin{align}\label{d}
    d_b=\sum_i e^i\wedge\nabla_{e_i}\,\,,\,\,\,\, \delta_b=-\sum_j i_{e_j}\nabla_{e_j}\,,\,\,\,i,j=1,\ldots, q.
\end{align}

\noindent For a basic $k$-form $\alpha$ let us define a normal $k$-tensor field $\rho$ on $(M,\mathcal{F})$ given at each point $x\in Q$ by
\begin{align}\label{rho}
    \rho_\alpha(X_1,\ldots,X_k)=\sum_{i=1}^{q}\sum_{j=1}^{k} (R(e_i,X_j)\alpha) (X_1,\ldots,X_{j-1},e_i,X_{j+1},\ldots,X_k),
\end{align}
where $X_j\in(\Gamma Q)_x$ and $R$ is the normal Riemann curvature operator. Using (\ref{d}) one can get the Weitzenb{\"o}ck formula for the Laplacian of a basic $k$-form $\alpha$ 
\begin{align}\label{Delta}
    \Delta\alpha=-\mathrm{tr}(\nabla^2\alpha)+\rho_\alpha,
\end{align}
for more details see \cite{po}.

The main result of A. Lichnerowicz in \cite{li} has its foliated counterpart, however, first we need to show that (\ref{Delta}) is equal to the Laplacian used in \cite{li}. A careful look at how the Riemann curvature operator acts on forms, one can from (\ref{rho}) obtain
\begin{flalign}
    \rho_\alpha(X_1,\ldots,X_k)=&-\sum_{i,j}\sum_l^{l\neq j}\alpha(X_1,\ldots,R(e_i,X_j)X_l,\ldots,e_i,\ldots,X_k)\\
    &-\sum_{i,j}\alpha(X_1,\ldots,R(e_i,X_j)e_i,\ldots,X_k).\nonumber
\end{flalign}

\noindent Now, if we apply the Ricci identity we get
\begin{flalign}\label{rho2}
    \rho_\alpha(X_1,\ldots,X_k)=&-\sum_{i,j,l}^{l\neq j}\sum_{m=1}^q R(e_i,X_j,e_m,X_l)\alpha(X_1,\ldots,e_m,\ldots,e_i,\ldots,X_k)\\
    &+\sum_{i,j}{\text{ Ric}}(e_i,X_j)\alpha(X_1,\ldots,e_i,\ldots,X_k),\nonumber
\end{flalign}
where in the second line we have used the fact that 
\begin{align*}
    R(e_i,X_j,e_m,e_i)=-R(e_i,X_j,e_i,e_m)=-{\text{Ric}(X_j,e_m)}=-{\text{Ric}(e_m,X_j)}.
\end{align*}

Let $\nabla$ be a torsionfree foliated connection in $(N(M,\mathcal{F}),\mathcal{F}_N)$, if one replaces $\rho_\alpha$ in (\ref{Delta}) with the one in (\ref{rho2}), it is clear to be the same as the Laplacian used by Lichnerowicz in \cite{li}, hence the basic Laplacian of a normal $p$-tensor field can be written locally using the Ricci tensor and Riemann curvature tensor of the foliated manifold 
\begin{align*}
    \Delta_b U_{\alpha_1\ldots\alpha_p}=-\nabla^i\nabla_i U_{\alpha_1\ldots\alpha_p}+\sum_{j}\sum_{A}R_{\alpha_j A}U_{\alpha_1\ldots^A\ldots\alpha_p}-\sum^{j\neq k}_{j,k}\sum_{A, B}R_{\alpha_jA,\alpha_kB}U_{\alpha_1\ldots^A\ldots^B\ldots\alpha_p}.
\end{align*}

\begin{theorem}\label{thm4}
Let $T$ be a normal tensor field. If the covariant derivative of the normal tensor field $T$ vanishes, then for any normal tensor field $U$ the following holds

\begin{align*}
    \Delta_b(T\otimes U)=T\otimes \Delta_b U.
\end{align*}
\begin{proof}\normalfont
Let $T$ and $U$ be two normal tensor fields with the ranks $p$ and $q$ respectively. Applying the basic Laplacian on the tensor product $T\otimes U$ and using $\nabla T=0$, we obtain 
\begin{align*}
    \Delta_b(T\otimes U)=T\otimes \Delta_b U+ V\otimes U + W.
\end{align*}

\noindent
where the tensors $V$ and $W$ are 
\begin{flalign*}\label{no1}
    V_{\mu_1\ldots\mu_p}=\sum_{j}\sum_{A}R_{\mu_j A}T_{\mu_1\ldots^A\ldots\mu_p}-\sum^{j\neq k}_{j,k}\sum_{A, B}R_{\mu_jA,\mu_kB}T_{\mu_1\ldots^A\ldots^B\ldots\mu_p},\\
    W_{\mu_1\ldots\mu_p\nu_1\ldots\nu_q}=-\sum_{j,k}\sum_{A, B}R_{\mu_jA,\mu_kB}T_{\mu_1\ldots^A\ldots\mu_p}U_{\nu_1\ldots^B\ldots\nu_q}.
\end{flalign*}
So the aim is now to show that $W$ and $V\otimes U$ vanish. Applying the Ricci identity on the normal tensor $T$ we get

\begin{align}
(\nabla_\alpha\nabla_\beta - \nabla_\beta\nabla_\alpha)T_{\mu_1\ldots\mu_p}=\sum_j\sum_{A} R_{{\mu_j}A,\alpha\beta}T_{\mu_1\ldots^A\ldots\mu_p}=0
\end{align}
which implies that $W$ vanishes. The basic Laplacian commutes with the contraction and if in (\ref{no1}) we contract $\mu_j$ with $\beta$ and change $\alpha$ to $\mu_k$ we obtain
\begin{align*}
    \sum_{A}R_{\mu_kA}T_{\mu_1\ldots^A\ldots\mu_p}=\sum^{j\neq k}_j \sum_{A, B} R_{\mu_kA,\mu_kB}T_{\mu_1\ldots^A\ldots^B\ldots\mu_p}=0
\end{align*}

\noindent
which results in $V=0$ and this completes the proof.
\end{proof}
\end{theorem}

\noindent As a consequence, we obtain the following corollary.

\begin{corollary} If $\mathcal{T}$ is a linear mapping of the module of normal $r$-tensor
fields into normal $s$-tensor fields defined by a tensor field $T$ with $\Delta_bT =
0$, then $\mathcal{T}$ commutes with $\Delta_b$.

\end{corollary}

Kraines noticed that the Chern decomposition theorem, see \cite{ch}, can
be applied in the context she studied. Thus a harmonic form $\omega$ can be
represented as
\begin{align*}
    \omega=\sum^{[p/4]}_{i=0}L^i\omega^{p-4i}_e,
\end{align*}
\noindent
and then the forms $L^i\omega^{p-4i}_e$ must be harmonic.\newline

The original Chern’s proof of the decomposition theorem is just a
very subtle linear algebra plus the Hodge decomposition theorem for
differential forms. The theory of harmonic basic forms for compact
Riemannian foliated manifolds permits us to extend the result to basic
forms, and we get a basic counterpart of Kraines’ Theorem
$3.6$. 

\begin{theorem}\label{thm5}
Let $(M,\mathcal{F})$ be a compact Riemannian foliated manifold of
codimension $4q$. If the foliation $\mathcal{F}$ is cohomologically taut and transversely quaternionic K{\"a}hler
then the basic Betti numbers $B^i_\mathcal{F}$
of $(M,\mathcal{F})$ satisfy the inequalities:

\begin{align*}
    B^i_\mathcal{F}\leq  B^{i+4}_\mathcal{F}\leq \ldots\leq  B^{i+4r}_\mathcal{F}
\end{align*}
\noindent
for $i + 4r \leq q + 1, i = 0, 1, 2$ or $3$.

\begin{proof}\normalfont
Since the basic $4$-form $\Omega$ is invariant under $G$, so are the subspaces $L^i\bigwedge^{p-4i}_e$ of $\bigwedge^p$. Therefore, they can be written as sum of the subspaces $W_1,\ldots,W_s$ and the projection of a harmonic form into these subspaces is again harmonic.
\end{proof}
\end{theorem}

\noindent The result above holds  for a compact $3$-Sasakian manifold $M^{4n+3}$, and moreover K. Galicki and S. Salamon formulated and demonstrated the following theorem, \cite{gal}.

\begin{theorem}
    Let $(M^{4n+3)},g)$ be a compact $3$-Sasakian manifold. Then the odd Betti numbers $\mathcal{B}_{2k+1}$ of $M$ are all zero for $0\leq k\leq n$.
\end{theorem}

 The previous results of the section, combined with the proof provided in \cite{fu} allow us to formulate the following foliated version of \cite[Theorem 3.22]{fu}.

\begin{theorem}\label{th5}
 Let $(M, g, Q, \mathcal{F})$ be a cohomologically taut quaternionic K{\"a}hler
foliated manifold of codimension $4q$. Then
\begin{enumerate}[1)]
    \item for any $k < q$ the linear map $L\colon H^k(M,\mathcal{F}) \rightarrow H^{k+4}(M,\mathcal{F})$ is
injective,
    \item and there is the direct sum decomposition
$H^k(M,\mathcal{F}) = \sum_{0\leq s\leq [k/4]}L^s H^{k-4s}_e(M,\mathcal{F})$, \break $k\leq q+3$.
\end{enumerate}
\end{theorem}

\begin{example}

In the paper \cite{fu} A. Fujiki presents a generalization of Kraines's Theorem, cf. \cite[Theorem 3.23]{fu}. The proof uses two important facts, namely that the total space of the twistor bundle is a compact K{\"a}hler manifold and that the pull-back mapping induces an injection of the cohomology of the base into the cohomology of the total space. This fact is due to the Leray-Hirsch theorem. It has been very tempting to formulate and prove a foliated version of Fujiki's theorem. A careful reading of Salomon's theorem, cf. \cite[Theorem 6.1]{sa}, confirms that the canonical foliation of the transverse twistor space of a transversely quaternionic foliation is transversely K\"ahler. One should also mention that Aziz El Kacimi has demonstrated a foliated version of the Calabi conjecture, cf. \cite{el}. However, there is no version of the Leray-Hirsch theorem for foliated bundles and their basic cohomology. On the other hand, it is possible for one interesting class of foliated manifolds, namely $3$-$(\alpha,\delta)$-Sasaki and in particular $3$-Sasakian manifolds (i.e. $\alpha=\delta=1$). Indeed, their canonical foliation is transversely quaternionic K\"ahler.

\begin{theorem}
 Let $M$ be a compact $3$-$(\alpha,\delta)$-Sasaki manifold of dimension $3+4q$ and $\alpha\delta>0$.  Then
\begin{enumerate}[(1)]
    \item for any $0 \leq k < q$ the induced linear map $L^{q-k}\colon H^{2k}(M,\mathcal{F}) \longrightarrow
H^{4q-2k}(M,\mathcal{F})$ is an isomorphism.
    \item for any $k \geq 0$ we have the direct sum decomposition
$H^{2k}(M,\mathcal{F}) = \sum_{r} L^r H^{2k-4r}_e (M, \mathcal{F}),\break (k-q) \leq r \leq [k/2]$.
\end{enumerate} 
\begin{proof}
    As $M$ is a compact $3$-$(\alpha,\delta)$-Sasaki manifold with $\alpha\delta>0$, from \cite{Ag20} we see that its characteristic foliation $\mathcal{F}$ has positive normal scalar curvature, then from the results of El Kacimi \cite{el} and Salamon \cite{sa} it follows that the transverse twistor space is transversely K\"ahler. As the transverse twistor space in this case is a trivial bundle, there is no need for a foliated version of Leray-Hirsch theorem and the rest follows from the proof of \cite[Theorem 3.23]{fu}. 
\end{proof}

\end{theorem}
Obviously, the same result holds for compact $3$-Sasakian manifolds as from \cite{bg1} we have that its characteristic foliation $\mathcal{F}$ has positive normal scalar curvature $2(2q+1)(4q+3)$ and the rest of the proof is similar.

\end{example}

\section{Orbifolds}\label{sec9}

In 1956, I. Satake \cite{sat} introduced a new generalization of the notion of manifolds that he named \break $V$-manifolds. Currently, due to W. Thurston \cite{th}, they are known as orbifolds and have applications in both mathematics and physics, especially in string theory. It is a well-known result that having a Riemannian foliation with compact leaves, its leaf space can be given a structure of an orbifold and any orbifold can be realized as the leaf space of a Riemannian foliation, cf. \cite{gi}. This enables us to reformulate some of our results for orbifolds. In this section, we follow the notation and borrow some notions from \cite{wo-1}, which can be consulted for more discussions on the subject.

Let $X$ be a topological space, $\tilde{U}\subset\mathbb{R}^n$ be a connected open subset, $\Gamma$ be a finite group of smooth diffeomorphisms of $\tilde{U}$, and $\phi\colon\tilde{U}\longrightarrow X$ be a map which is $\Gamma$-invariant and induces a homeomorphism of $\tilde{U}/\Gamma$ onto an open subset $U\subset X$. The triple $(\tilde{U},\Gamma,\phi)$ is called an $n$-dimensional orbifold chart on X.

An embedding $\lambda\colon(\tilde{U},\Gamma,\phi)\longrightarrow (\tilde{V},\Delta,\psi)$ between two orbifold charts is a smooth embedding $\lambda:\tilde{U}\longrightarrow \tilde{V}$ which satisfies $\psi\circ\lambda=\phi$.

Let $\mathcal{A}=\{(\tilde{U}_i,\Gamma_i,\phi_i)\}_{i\in I}$ be a family of such charts, it is called an orbifold atlas on $X$, if it covers $X$ and any two charts are locally compatible in the following sense:
given two charts $\{(\tilde{U}_i,\Gamma_i,\phi_i)\}_{i=1,2}$ and $x\in U_1\cup U_2$, there exists an open neighborhood $U_3\subset U_1\cup U_2$ containing $x$ and a chart $(\tilde{U}_3,\Gamma_3,\phi_3)$, $U_3=\phi_3(\tilde{U}_3)\subset X$ such that it can be embedded into the other two charts. As in the case of manifolds, one can define a maximal atlas.

\begin{definition}
A Hausdorf paracompact topological space $X$ together with a maximal orbifold atlas $\mathcal{A}$ is called a smooth $n$-dimensional orbifold.
\end{definition}

Let $X$ be the orbifold associated to the foliated Riemannian manifold $(M,\mathcal{F})$ with compact leaves. The foliated geometrical structures on $M$ are in one-to-one correspondence with the orbifold geometrical structures on $X$, e.g., any bundle-like Riemannian metric of $M$ induces an orbifold Riemannian metric on $X$ and vice versa, cf. \cite{wo-1}.

\begin{definition}
Let $F$ be a smooth manifold. An orbifold $E$ is called an orbifold frame bundle over the orbifold $X$ with standard fiber $F$ if

\begin{itemize}
    \item[i)] there exists a smooth orbifold map $p\colon E\longrightarrow X$,  
    \item[ii)] there exists an orbifold atlas $\mathcal{A}=\{(\tilde{U}_i,\Gamma_i,\phi_i)\}$ of $X$,
    \item[iii)] let $V_i=p^{-1}(U_i)$ and $\tilde{V}_i=\tilde{U}_i\times F$, then there exist a group $\Lambda_i$ of fiber preserving diffeomorphisms of $\tilde{V}_i$ and a homeomorphism $\psi_i\colon \tilde{V}_i/\Lambda_i\longrightarrow V_i$ such that $\{(\tilde{V}_i,\Lambda_i,\psi_i)\}$ form an atlas of the orbifold $E$
    \item[iv)] and the following diagram is commutative
    
    \begin{displaymath}
  \begin{tikzcd}
   \tilde{U}_i\times F \arrow[rr,  "\tilde{p}=p\times id "] \arrow[d,]
   && \tilde{U}_i \arrow[d, ]\\
   \tilde{V}_i/\Lambda_i \arrow[d, "\phi_i "]
   && \tilde{U}_i/\Gamma_i \arrow[d, "\psi_i"]\\
   V_i \arrow[rr, "p"]
   && U_i
 \end{tikzcd}
 \end{displaymath}
\end{itemize}
\end{definition}

The tangent bundle of an orbifold can be constructed as follows. Let $\{(\tilde{U}_i,\Gamma_i,\phi_i)\}$ be the orbifold atlas on $X$. Take $\tilde{V}_i=T\tilde{U}_i=\tilde{U}_i\times\mathbb{R}^n$ and for the group $\Sigma_i$ of local transformations take $\Sigma_i=\{d\gamma_i:\gamma\in \Gamma_i\}$ and the quotient map $\tilde{V}_i\longrightarrow \tilde{V}_i/\Sigma_i$ can be taken as $\psi_i$. It can be shown that $\{(\tilde{V}_i,\Sigma_i,\psi_i)\}$ is an orbifold atlas for $TX$.

Denote by $L(X)$ the linear frame bundle of the orbifold $X$. It is constructed similarly as the tangent bundle $TX$ with fiber $F$ now being $GL(n)$ instead of $\mathbb{R}^n$. It is well-known that $L(X)$ is in fact a manifold and many of the geometrical structures on $X$ can be realized as its reduction, e.g., having an orbifold Riemannian metric on $X$ is equivalent to a choice of an orbifold $O(n)$ reduction of $L(X)$. In particular an orbifold quaternionic K{\"a}hler structure on $X$ can be introduced by an orbifold $Sp(n).Sp(1)$ reduction of $L(X)$. One can also construct a dual vector bundle of any orbifold vector bundle, more details can be found in another article written by the second author \cite{wo-1}. Moreover, the tensor product, skewsymmetric product and exterior product of orbifold vector bundles over a given orbifold  can be defined and an orbifold differential $k$-form on the orbifold $X$ can be taken as the section of $\bigwedge^k T^*X$, i.e., the $k$-th exterior power of the cotangent bundle of $X$. 

Let $\Omega^k(X)$ be the space of all orbifold differential $k$-forms on the orbifold $X$, having the differential $d\colon\Omega^k(X)\longrightarrow \Omega^{k+1}(X)$, one can define the orbifold de Rham cohomology group $H^*_{DR}(\Omega^*,d)$ similarly as for manifolds, which for the sake of brevity from now on we shall use $H^*_{DR}(X)$, for more details, see \cite{alr}.

\begin{theorem}\cite[Theorem 2.15]{mm}
Let $\mathcal{F}$ be a foliation of codimension $q$ of a manifold $M$ such that any leaf of $\mathcal{F}$ is compact with finite holonomy group. Then the space of leaves $M/\mathcal{F}$ has a canonical structure of dimension $q$. The isotropy group of a point in $M/\mathcal{F}$ is its holonomy group.
\end{theorem}

 An effective orbifold differential $k$-form can be defined similarly as for manifolds and we denote its orbifold cohomology group by $H^*_{e,DR}(X)$. It can be shown that if $X$ is the orbifold associated to the foliated Riemannian manifold $M$, the de Rham orbifold cohomology group on $X$ is actually isomorphic to the basic cohomology group of $M$, hence, we can reformulate Theorem \ref{th5}.

\begin{theorem}
Let $(M, g, Q, \mathcal{F})$ be a cohomologically taut quaternionic K{\"a}hler
foliated manifold of codimension $4q$ and let $X$ be its associated orbifold with its de Rham cohomology group denoted by $H^*_{DR}(X)$. Then
\begin{enumerate}[(1)]
    \item for any $k < q$ the linear map $L \colon H^k_{DR}(X) \rightarrow H^{k+4}_{DR}(X)$ is
injective,
    \item and there is the direct sum decomposition
$H^k_{DR}(X) = \sum_{0\leq s\leq [k/4]}L^s H^{k-4s}_{e,DR}(X)$,\,\,\,\,\,  $k\leq q+3$.
\end{enumerate}
\end{theorem}

Let $M^{4n+3}$ be a $3$-Sasakian manifold, and let $\mathcal{F}_\xi$ denote its $3$-dimensional quasi-regular foliation. This is perhaps the most important and natural example for this section. Since $\mathcal{F}_\xi$ is quasi-regular, the leaf space $\mathcal{O}:=M/\mathcal{F_\xi}$ is an orbifold with quaternionic K\"ahler structure, cf. \cite{bg1}, while its regular case counterpart was proved almost 20 years earlier, as we mentioned in Subsection \ref{sub2.3}. Now following a similar argument as in Appendix B, the projection $\pi\colon M\rightarrow \mathcal{O}$, induces an isomorphism $H^*_{DR}(\mathcal{O})\cong H^*(M,\mathcal{F})$. The relation between the cohomology of the quasi-regular compact $3$-Sasakian manifold $M$ and the orbifold $\mathcal{O}$, gives also a relation between the Betti numbers.

\begin{prop}\cite[Proposition 13.5.5]{bg2}
    Let $(M^{4n+3},g)$ be a compact $3$-Sasakian manifold and $\mathcal{O}=M/\mathcal{F_\xi}$. Then  $\mathcal{B}_{2p}(M)=\mathcal{B}_{2p}(\mathcal{O})-\mathcal{B}_{2p-4}(\mathcal{O})$ for $p\leq 2n+1.$
\end{prop}

\noindent An important observation is that a harmonic $k$-form
on $M$ is the lift of a primitive harmonic form on $\mathcal{O}$,
$0 \leq k \leq 2n + 1$, cf. \cite{bg2}.

A classical result in differential geometry, proven in \cite{tho}, asserts that a compact four-dimensional self-dual manifold with positive scalar curvature is either isometric to $S^4$ or $\mathbb{CP^2}$ with the Fubini-Study metric. In \cite{Hi}, N. Hitchin constructed a family of self-dual Einstein orbifold metrics in dimension four, parametrized by an integer $k\geq 3$, which from now on we shall denote by $\mathcal{O}_k$. One of their main properties is that they are invariant under an $SO(3)$ action. Now consider their Konishi bundle, which is an $SO(3)$ orbifold principal bundle with a $3$-Sasakian structure, and we shall denote it by $S_k$; then by a result of K. Grove et. al. \cite{gr}, $S_k$ is a smooth $3$-Sasakian manifold for each $k$. Therefore, we have a very nice example for this section, the most trivial example being the Hopf fibration $S^7\rightarrow S^4$.

\appendix
\setcounter{equation}{0}
\renewcommand{\theequation}{A.\arabic{equation}}
\section{Quaternions and differential forms}\label{App1}

Let $\mathbb{H}^n$ be the $n$-dimensional right module over the field of quaternions $\mathbb{H}$. The canonical bilinear form on $\mathbb{H}^n$ is defined as
\begin{equation*}
    \langle u,v\rangle=\frac{1}{2}\sum^{n}_{i=1}\left(u_i\Bar{v}_i+v_i\Bar{u}_i\right)
\end{equation*}
where $u=\left(u_1,\ldots,u_n\right), v=\left(v_1,\ldots,v_n\right)\in \mathbb{H}^n$. The $2$-form $\langle .\,,.\,\rangle$ is a scalar product on the $4n$-dimensional real vector space $\mathbb{H}^n$. The group $Sp(n)$ can be defined as the linear group preserving the ``symplectic product" $\left(u,v\right)=\sum^{n}_{i=1}u_i\Bar{v}_i$ on $\mathbb{H}^n$. Then $\langle u,v\rangle=\frac{1}{2}\{\left(u,v\right)+\left(v,u\right)\}$ thus immediately the natural action of $Sp(n)$ preserves the scalar product $\langle .\,,.\,\rangle$. The group $Sp(1)$ is identified with the quaternions of length $1$, so its right action preserves the product $\langle .\,,.\,\rangle$.

As $\mathbb{H}$ can be identified with $\mathbb{R}^4$, it makes it possible to write any quaternion $x\in \mathbb{H}$ as \break $x=x^0\textbf{1}+x^1\textbf{i}+x^2\textbf{j}+x^3\textbf{k}$ where $\textbf{1},\textbf{i},\textbf{j},\textbf{k}$ form the standard base of $\mathbb{H}$ as a real vector space. The right multiplication by $\textbf{i}, \textbf{j}$ and $\textbf{k}$ define three complex structures on $\mathbb{H}^n$ denoted by the same letters, respectively. In turn, they permit us to define three skew-symmetric $2$-forms
\begin{equation*}
    \Omega_I\left(u,v\right)=\langle u\textbf{i},v\rangle,\,\,\,\,\Omega_J\left(u,v\right)=\langle u\textbf{j},v\rangle,\,\,\,\,\Omega_K\left(u,v\right)=\langle u\textbf{k},v\rangle.
\end{equation*}

\noindent In \cite{kr} V.Y. Kraines demonstrates that the $4$-form $\Omega$ on $\mathbb{H}^n$ defined by
\begin{equation}\label{omegaH}
    \Omega=\Omega_I\wedge \Omega_I+\Omega_J\wedge \Omega_J+\Omega_K\wedge \Omega_K.
\end{equation}
is invariant under the natural action of $Sp(n).Sp(1)$.

Let $(\mathbb{H}^n)^\prime$ be the dual space of $\mathbb{H}^n$ as the quaternionic vector space. Let $z_1,\ldots,z_n$ be a basis of $(\mathbb{H}^n)^\prime$. Each $z_\alpha$ can be represented as

\begin{equation*}
    z_\alpha=a_\alpha\textbf{1}+b_\alpha\textbf{i}+c_\alpha\textbf{j}+d_\alpha\textbf{k},
\end{equation*}
thus the $1$-forms $a_1,\ldots,a_n,b_1,\ldots,b_n,c_1,\ldots,c_n,d_1,\ldots,d_n$ form a basis of $(\mathbb{H}^n)^\prime$ as a real vector space. Then, $\Omega=\sum_{\alpha,\beta,\gamma,\delta}a_\alpha\wedge b_\beta\wedge c_\gamma\wedge d_\delta,$ hence $\Omega^n\neq 0.$ It follows that the elements of $\bigwedge^* \left(\mathbb{H}^n\right)^\prime$ can be written as linear combination of simple $p$-forms $\omega=\beta_1\wedge\ldots\wedge\beta_p$, where the $1$-forms $\beta_i$ can be one of the $\{a_\alpha, b_\alpha, c_\alpha, d_\alpha\}$ basis. On $\bigwedge^* \left(\mathbb{H}^n\right)^\prime$ we can define $3$ operators $*, L$ and $\Lambda$ as follows, cf. \cite{kr}:
If $\omega$ is a simple $p$-form, then $*\omega$ is the simple $(4n-p)$-form that

\begin{equation*}
    \omega\wedge *\omega=a_1\wedge b_1\wedge c_1\wedge d_1\wedge \ldots \wedge a_n\wedge b_n\wedge c_n \wedge d_n.
\end{equation*}
Then $*$ is extended by $\mathbb{R}$-linearity to $\bigwedge \left(\mathbb{H}^n\right)^\prime$. Moreover, $**\omega=\omega$. For an arbitrary form $\omega$ we define $L$ and $\Lambda$ operators as $  L\omega=\Omega\wedge\omega\; \text{and}\; \Lambda\omega=*(\Omega\wedge*\omega).$ Thus $L$ and $\Lambda$ are linear transformations
\begin{center}
  $L\colon\bigwedge^p \left(\mathbb{H}^n\right)^\prime \rightarrow \bigwedge^{p+4} \left(\mathbb{H}^n\right)^\prime,$\\
    $\Lambda\colon\bigwedge^p \left(\mathbb{H}^n\right)^\prime \rightarrow \bigwedge^{p-4} \left(\mathbb{H}^n\right)^\prime.$
   \end{center} 
    
\noindent On $\bigwedge^p \left(\mathbb{H}^n\right)^\prime$ we define a bilinear form by
    \begin{equation*}
        \left(\omega,\omega^\prime\right)=*\left(\omega\wedge*\omega^\prime\right).
    \end{equation*}
    where $\omega,\omega^\prime\in\bigwedge^p(\mathbb{H}^n)^\prime$. Then
    \begin{equation*}
        \left(L\omega,\omega^\prime\right)=\left(\omega,\Lambda\omega^\prime\right).
    \end{equation*}
    for any $\omega\in \bigwedge^p \left(\mathbb{H}^n\right)^\prime$ and $\omega^\prime\in \bigwedge^{p+4} \left(\mathbb{H}^n\right)^\prime$.\newline
    
\noindent Kraines shows that $L\colon\bigwedge^p \left(\mathbb{H}^n\right)^\prime \longrightarrow \bigwedge^{p+4} \left(\mathbb{H}^n\right)^\prime$ is an isomorphism into for $p+4\leq n+1$. Next, she defines effective forms: a $p$-form $\omega$ is called effective if $\Lambda\omega=0$. The space of all effective $p$-forms is denoted by $\Lambda^p_e$. With all these notions in place, she formulates and proves the following decomposition, cf. \cite[Theorem 2.6]{kr}.

\begin{theorem}\label{Thm1}
There is the following direct sum decomposition of $\bigwedge^p \left(\mathbb{H}^n\right)^\prime$ for $p\leq n+1$, $r=[p/4]$,

\begin{center}
 $\bigwedge^{p} \left(\mathbb{H}^n\right)^\prime=\Lambda^p_e+L\Lambda^{p-4}_e+\ldots+L^r\Lambda^{p-4r}_e.$   
\end{center}

\end{theorem}

\setcounter{equation}{0}
\renewcommand{\theequation}{B.\arabic{equation}}

\section{Hodge theory for basic forms}\label{sec5}

Let $(M^m, \mathcal{F}, g)$ be a foliated Riemannian manifold, $m=p+4n$, where the dimension and codimension of the foliation are $p$ and $q=4n$ respectively. On a foliated Riemannian manifold, apart from the complex of differential $k$-forms $A^k(M)$, there exists a complex of basic $k$-forms. The former is quite well known; for the reader's convenience, we recall the latter and discuss the definitions and properties of Hodge theory for basic forms. On a foliated Riemannian manifold $(M,g,\mathcal{F})$ the set of all basic k-forms is 
\begin{align*}
A^k(M,\mathcal{F}) = \{\alpha \in A^k(M) : i_X\alpha =0,\, i_Xd\alpha = 0\, \text{for all vectors}\, X \in T\mathcal{F}\}
\end{align*}
which is a subcomplex of $A^k(M)$ and we denote its cohomology by $H^k(M,\mathcal{F})$. The $k$-th basic Betti number is then defined as $B^k_\mathcal{F}=dim H^k(M,\mathcal{F})$. Denote by $\pi:M\rightarrow M/\mathcal{F}$, the canonical projection map, where $M/\mathcal{F}$ is the leaf space. Let $(U,\phi)$ be a foliated chart with coordinates $\phi(x_1,\ldots,x_p,y_1,\ldots,y_q)$, then a basic $k$-form $\alpha$ can be written as
\begin{align*}
    \alpha=\sum_{1\leq i_1<\cdots<i_k\leq q}\alpha_{i_1\ldots i_k(y_1,\ldots,y_q)dy_{i_1}\wedge\ldots\wedge dy_{i_q}}.
\end{align*}

Let $A^*(M/\mathcal{F})$ denote the complex of differential forms on $M/\mathcal{F}$, then the pullback of the projection map $\pi^*:A(M/\mathcal{F})\rightarrow A^*(M,\mathcal{F})$ is an isomorphism of differential complexes. Moreover, we have the isomorphism $H^*(M/\mathcal{F})\cong H^*(M,\mathcal{F})$.

The restriction of the bundle-like metric to the normal bundle of the foliation of the Riemannian foliated manifold $(M, g, \mathcal{F})$ defines $\Bar{*}$ operator, cf. \cite{to},
\begin{align*}
    \Bar{*}\colon A^k(M,\mathcal{F})\rightarrow A^{q-k}(M,\mathcal{F}).
\end{align*}

\noindent On the Riemannian manifold $(M, g)$ we have the $*$-operator acting
on the complex of smooth forms:
\begin{align*}
    *\colon A^k(M)\rightarrow A^{m-k}(M)
\end{align*}
On the subcomplex $A^k(M,\mathcal{F})$ of basic forms these two operators are
related by the formula
\begin{align*}
    \Bar{*}\alpha=(-1)^{p(q-k)}*(\alpha\wedge\chi_\mathcal{F}),
\end{align*}
for any $\alpha \in A^k(M,\mathcal{F})$, where $\chi_\mathcal{F}$ is the volume form of leaves.\newline

\noindent In $A^k(M,\mathcal{F})$ we have the standard scalar product
\begin{align*}
    \langle \alpha,\beta\rangle_b=\int_M \alpha\wedge \Bar{*}\beta\wedge\chi_\mathcal{F},
\end{align*}
which is the restriction of the standard scalar product on $A^k(M)$.
A Riemannian foliation on a compact manifold is said to be taut if
there exists a Riemannian metric that makes all its leaves minimal
submanifolds. Tautness is characterized by the nonvanishing of the top
dimensional basic cohomology, i.e., $H^q(M,\mathcal{F}) \neq 0$. In this case we say that the foliation is cohomologically taut. In fact, this Riemannian metric can be chosen to be bundle-like. Moreover, one can
make the modification only along leaves, cf. \cite[Chapter $7$]{to}.

The formal adjoint $\delta_b$ of $d$ in the complex $A^k(M,\mathcal{F})$ with the scalar product $\langle .\,,.\,\rangle_b$ is the operator
\begin{align*}
    \delta_b=(d-\kappa\wedge)^{\Bar{*}}\colon A^k(M,\mathcal{F})\rightarrow A^{k-1}(M,\mathcal{F}),
\end{align*}
where $\kappa$ is the mean curvature form of the leaves, and
\begin{align*}
 (d-\kappa\wedge)^{\Bar{*}}(\beta) = (-1)^{q(k+1)+1}\Bar{*}(d-\kappa)\Bar{*}\beta,   
\end{align*}
for any $\beta \in A^k(M,\mathcal{F})$. If the leaves of $\mathcal{F}$ are minimal submanifolds for the bundle-like metric $g$, then $\kappa = 0$ and $\delta_b = d^{\Bar{*}}$. We define the basic Laplacian as
\begin{align*}
    \Delta_b=\delta_bd+d\delta_b.
\end{align*}

\noindent A basic form $\alpha$ is called harmonic iff $\Delta_b \alpha = 0$. Moreover, for compact Riemannian foliated manifolds asserts that $\alpha$ is
harmonic iff $d\alpha = 0 = \delta_b\alpha$.

\setcounter{equation}{0}
\renewcommand{\theequation}{C.\arabic{equation}}

\section{Foliated \texorpdfstring{$G$}{G}-structures}\label{App3}

The second author introduced the structure tensor of a transverse $G$-structure on a foliated manifold in \cite{wo-b}. We will first recall some of the definitions and results from \cite{wo-b} and then discuss the case of the Riemannian foliated manifold. 

Let $(M,\mathcal{F})$ be a foliated manifold of codimension $q$ and $N(M,\mathcal{F})$ be the normal bundle of the foliation. The bundle of linear frames of $N(M, \mathcal{F})$ is a $GL(q)$-principal bundle and we denote it by $L(M, \mathcal{F})$. Let $L$ be the total space of the bundle $L(M, \mathcal{F})$, it is also foliated with leaves that project to the leaves of $\mathcal{F}$ and we denote the foliation by $\mathcal{F}_L$. Let $G$ be a closed subgroup of $GL(q)$

\begin{definition}\cite{wo-b}
    \begin{itemize}
        \item A \emph{transverse $G$-structure $B(M,G;\mathcal{F})$} is a $G$-reduction of the principal bundle $L(M,\mathcal{F})$.
        \item A \emph{foliated G-structure $B(M,G;\mathcal{F})$} is a $G$-reduction of $L(M,\mathcal{F})$ whose total space $B$ is $\mathcal{F}_L$-saturated (i.e. if it contains a point then it also contains the whole leaf of $\mathcal{F}_L$ passing through the point).
    \end{itemize}
\end{definition}

\begin{definition}\cite{wo-b}
    \begin{enumerate}
        \item A connection in a transverse $G$-structure is called a transverse connection.
        \item A connection in a foliated $G$-structure is called:
        \begin{itemize}
            \item basic if its connection form vanishes on vectors tangent to the foliation $\mathcal{F}_L$,
            \item foliated if its connection form is basic.
        \end{itemize}
    \end{enumerate}
\end{definition}

\begin{lemma}\cite[Lemma 2]{wo-b}
    A foliated $G$-structure $B(M, G; \mathcal{F})$ admits a foliated $G$-connection if and only if the corresponding $G$-structure on the transverse manifold admits an $\mathcal{H}_\mathcal{U}$-invariant\break $G$-connection.
\end{lemma}

Let $\pi\colon L\rightarrow M$ denote the projection from the foliated frame bundle $L(M, \mathcal{F}_L)$ to the foliated manifold $(M,\mathcal{F})$ and $\pi_N: TM\rightarrow N(M,\mathcal{F)}$ be the natural projection from the tangent bundle of the foliated manifold to the normal bundle of the foliation. We define an $\mathbb{R}^q$-valued canonical $1$-form on the total space $L$ of $L(M, \mathcal{F}_L)$ or $B$ of $B(M, G;\mathcal{F})$ as follows
\begin{align}
    \theta_p=p^{-1}\pi_Nd_p\pi,\,\,\,x=\pi(p)
\end{align}
where $p^{-1}\colon N_x(M,\mathcal{F})\rightarrow \mathbb{R}^q$  is the inverse of the natural isomorphism defined by the linear normal frame $p$.\newline

Now for a foliated manifold $(M,\mathcal{F})$ with a transverse $G$-structure $B(M,G;\mathcal{F})$, fix a subbundle $Q$ of $TM$ such that it is supplementary to $T\mathcal{F}$. A $q$-dimensional subspace $H_p$ of the tangent space $T_p B$ is called horizontal if $d\pi(H_p)=Q_{\pi(p)}$. The following map can be defined for any horizontal subspace $H$
\begin{align}
    &C_H\colon\mathbb{R}^q\wedge\mathbb{R}^q\rightarrow \mathbb{R}^q,\\\nonumber
    &C_H(u\wedge v)=\langle X\wedge Y,d\theta\rangle.
\end{align}
where $X,Y\in H$ and $\theta(X)=u,\,\,\, \theta(Y)=v$. Let $H^\prime$ be another horizontal subspace with its corresponding map $C_{H^\prime}$. The maps $C_{H}$ and $C_{H^\prime}$ at any point $p\in B$ differ by an element of $\partial \mathsf{Hom}(\mathbb{R}^q,\mathfrak{g})$, where $\mathfrak{g}$ is the Lie algebra of $G$. Therefore for any point $p\in B$, the maps $C_H$ define a unique class in $\mathsf{Hom}(\mathbb{R}^q\wedge\mathbb{R}^q,\mathbb{R}^q)/\partial \mathsf{Hom}(\mathbb{R}^q,\mathfrak{g})$ denoted by $c(p)$. The tensor given in this way is called the structure tensor of the transverse $G$-structure $B(M,G;\mathcal{F})$ and is denoted by $c$.

\begin{prop}\cite[Proposition 5]{wo-b}
    Let $B(M,G;\mathcal{F})$ be a foliated $G$-structure on a foliated manifold $(M,\mathcal{F})$ and the first prolongation of the Lie algebra $\mathfrak{g}$ be trivial. Then if the structure tensor of $B(M, G; \mathcal{F})$ vanishes, there exists a torsion-free foliated connection in $B(M, G; \mathcal{F})$.
\end{prop}

Note that a foliated $G$-structure with $G$ being $O(q)$ is exactly a foliated Riemannian structure, i.e. we have a Riemannian foliated manifold $(M, \mathcal{F})$, and in this case, assumptions of the previous proposition are automatically satisfied and the torsion-free connection thus obtained is exactly the transverse Riemannian connection. Similarly, for Lie subgroups of $O(q)$, e.g. $Sp(q,\mathbb{R})$, with  $q=2n$, or $Sp(q).Sp(1)$ with $q=4n$, one gets transversely symplectic and transversely quaternion K\"ahler structures, respectively.

\end{document}